\documentclass[hypertex,11pt]{article}
\usepackage[pagebackref]{hyperref} 
\usepackage[backrefs, alphabetic,initials]{amsrefs}
\usepackage{amsfonts,amsmath,amscd,amssymb,amsthm}
\def\0{{\bar 0}}
\def\1{{\bar 1}}

\def\u{\overline{u}}
\def\w{\overline{w}}

\def\Ker {{\operatorname{Ker}\;}}
\def\End {{\operatorname{End}}}

\newcommand{\noi}{\noindent}
\newcommand{\ga}{\alpha}
\newcommand{\gb}{\beta}
\newcommand{\gc}{\gamma}

\newcommand{\gs}{\sigma}

\newcommand{\go}{\omega}
\newcommand{\gt}{\tau}

\newcommand{\gl}{\lambda}

\newcommand{\gk}{\kappa}
\newcommand{\gep}{\epsilon}

\newcommand{\fsl}{\mathfrak{sl}}

\newcommand{\fb}{\mathfrak{b}}

\newfont{\eufm}{eufm10 scaled\magstep1}

\newcommand{\bco}{\begin{conjecture}}
\newcommand{\ba}{\begin{alg}}
\newcommand{\ea}{\end{alg}}
\newcommand{\eco}{\end{conjecture}}
\newcommand{\bpf}{\begin{proof}}
\newcommand{\epf}{\end{proof}}
\newcommand{\bt}{\begin{theorem}}
\newcommand{\et}{\end{theorem}}

\newcommand{\br}{\begin{rem}}
\newcommand{\er}{\end{rem}}
\newcommand{\brs}{\begin{rems}}
\newcommand{\ers}{\end{rems}}
\newcommand{\bi}{\begin{itemize}}
\newcommand{\bl}{\begin{lemma}}
\newcommand{\bsul}{\begin{sublemma}}
\newcommand{\esul}{\end{sublemma}}
\newcommand{\bp}{\begin{proposition}}
\newcommand{\be}{\begin{equation}}
\newcommand{\bc}{\begin{corollary}}
\newcommand{\bexa}{\begin{example}}
\newcommand{\eexa}{\end{example}}
\newcommand{\bex}{\begin{exercise}}
\newcommand{\eex}{\end{exercise}}
\newcommand{\btab}{\begin{tab}}
\newcommand{\etab}{\end{tab}}
\newcommand{\ei}{\end{itemize}}
\newcommand{\el}{\end{lemma}}
\newcommand{\ep}{\end{proposition}}
\newcommand{\ee}{\end{equation}}
\newcommand{\ec}{\end{corollary}}
\newcommand{\Bc}{\begin{center}}
\newcommand{\Ec}{\end{center}}
\newcommand{\bh}{\begin{hyp}}
\newcommand{\eh}{\end{hyp}}
\newcommand{\bhs}{\begin{hyps}}
\newcommand{\ehs}{\end{hyps}}

\numberwithin{equation}{section}%
\newcommand{\lra}{\longrightarrow}

\begin{document}
\title{Table of Contents}
\newcommand{\rh}{\epsilon^{\widehat{\rho}}}
\newcommand{\pii}{\prod_{n=1}^\infty}
\newcommand{\sii}{\sum_{n= - \infty}^\infty}
\newtheorem{theorem}{Theorem}[section]
 \newtheorem{lemma}[theorem]{Lemma}
 \newtheorem{example}[theorem]{Example}
  \newtheorem{sublemma}[theorem]{Sublemma}
   \newtheorem{corollary}[theorem]{Corollary}
 \newtheorem{conjecture}[theorem]{Conjecture}
 \newtheorem{hyps}[theorem]{Hypotheses}
  \newtheorem{hyp}[theorem]{Hypothesis}
 \newtheorem{alg}[theorem]{Algorithm}
  \newtheorem{sico}[theorem]{Sign Convention}
 \newtheorem{tab}[theorem]{\quad \quad \quad \quad \quad \quad \quad \quad \quad \quad \quad \quad \quad Table }
  \newtheorem{rem}[theorem]{Remark}
  \newtheorem{rems}[theorem]{Remarks}
 \newtheorem{proposition}[theorem]{Proposition}
\font\twelveeufm=eufm10 scaled\magstep1 \font\teneufm=eufm10
\font\nineeufm=eufm9 \font\eighteufm=eufm8 \font\seveneufm=eufm7
\font\sixeufm=eufm6 \font\fiveeufm=eufm5
\newfam\eufmfam
\textfont\eufmfam=\twelveeufm \scriptfont\eufmfam=\nineeufm
\scriptscriptfont\eufmfam=\sixeufm \textfont\eufmfam=\teneufm
\scriptfont\eufmfam=\seveneufm \scriptscriptfont\eufmfam\fiveeufm
\newtheorem{exercise}{}
\setcounter{exercise}{0} \numberwithin{exercise}{section}
\newenvironment{emphit}{\begin{theorem}}{\end{theorem}}
\setcounter{Thm}{0} \numberwithin{Thm}{section} 
\title{Finitely generated, non-artinian monolithic modules.}
\author{Ian M. Musson\\Department of Mathematical Sciences\\
University of Wisconsin-Milwaukee\\ email: {\tt
musson@uwm.edu}}
\maketitle
\begin{abstract}
We survey noetherian rings $A$ over which the injective hull of every simple module is locally artinian.  Then we give a general construction for algebras $A$ that do not have this property. In characteristic 0, we also complete the classification of down-up algebras with this property which was begun in \cite{CarvalhoLompPusat} and \cite{CM2}.  \end{abstract}
\section  {\bf{\Large Introduction}}\label{Intro}
A module $M$ is {\it monolithic} if the intersection of all nonzero submodules of $M$ is nonzero.  The intersection of all nonzero submodules of  a monolithic module $M$ is a simple submodule known as the {\it lith} of $M.$  Thus monolithic modules have a unique lith!  This terminology is due to Roseblade \cite{R}, \cite{R1}.  It was pointed out to me by Ken Goodearl that monolithic modules are also known as subidirectly irreducible modules.   We consider the following property of
a noetherian ring  $A$.
$$(\diamond) \quad \mbox{ Every finitely generated monolithic } A\mbox{-module is artinian.} $$
\noi Equivalently, the injective hull of every simple $A$-module is  locally artinian.  Some history concerning property $(\diamond)$ is given in the introduction to \cite{CM2}. The property is not well understood, as is shown by the following quite baffling lists  of examples.\\ \\
The following rings $A$ have property $(\diamond)$.
\begin{itemize}
\item[{(A.0)}]  Commutative noetherian rings, and more generally PI and FBN rings \cite{J1}.
\end{itemize}
The next two examples are in fact PI rings.
\begin{itemize}
\item[{(A.1)}] The coordinate ring of the quantum plane, that is the algebra generated by elements $a, b$ subject to the relation $ab = q ba$ when $q \in K$ is a root of unity.
\item[{(A.2)}] The quantized Weyl algebra, that is the algebra generated by elements $a, b$ subject to the relation $ab - q ba = 1$ when $q \in K$ is a root of unity.
\item[{(A.3)}] The enveloping algebra $U(\fsl(2,K))$ where $K$ a field of characteristic 0, \cite{D}.
\item[{(A.4)}] The group rings $\mathbb{Z}G$ and $KG$ where $K$ is a field which is algebraic over  a finite field and $G$ is polycyclic-by-finite, \cite{J}, \cite{R1}.
\item[{(A.5)}] Prime noetherian rings of Krull dimension 1, \cite{CarvalhoLompPusat}, \cite{M1}.
\item[{(A.6)}]  There are simple noetherian, non-artinian rings for which any simple module is injective, and obviously these rings have property $(\diamond)$ \cite{CO}.
\end{itemize}
\vspace{0.4cm}
The following rings $A$ do not have property $(\diamond)$.
\begin{itemize}
\item[{(B.1)}] The coordinate ring of the quantum plane when $q \in K \backslash \{0\}$ is not a root of unity, \cite{CM2}.
\item[{(B.2)}] The quantized Weyl algebra, when $q \in K \backslash \{0\}$ is not a root of unity, \cite{CM2}.
\item[{(B.3)}] The enveloping algebra $U(\fb)$ over  an algebraically closed field of characteristic 0, when $\fb$ is finite dimensional,  solvable and non-nilpotent, \cite{CH}, \cite{M2}.
\item[{(B.4)}] The group algebra $KG$ where $K$ is a field which is not algebraic over  a finite field and $G$ is polycyclic-by-finite which is not nilpotent-by-finite, \cite{M1}.
\item[{(B.5)}] The Goodearl-Schofield example: a certain non-prime noetherian ring of Krull dimension 1, \cite{GS}.
\end{itemize}
\noi What has been lacking up to now is a general construction for finitely generated, non-artinian, monolithic modules.  In the next section we give such a construction under fairly mild conditions on $A$.  We show that examples (B.1)-(B.3) satisfy these conditions. We also  apply our construction to down-up
algebras in characteristic 0. Some open problems are given in the last section.\\ \\
I thank Allen Bell and Paula Carvalho for useful comments, and Toby Stafford for encouraging me to finish this paper.
\section{The construction.}
Let $K$ be a field. We make the following assumptions.
\begin{itemize}
\item[{(1)}] $A$  is a noetherian $K$-algebra without zero divisors.
\item[{(2)}] $w$ is a normal element of $A$.
\item[{(3)}] $J$ is a maximal left ideal such that $w-\mu \in J$ for some non-zero $\mu \in K$.
\end{itemize}
From (1) and (2) it follows that there is an automorphism $\gs$ of $A$ such that for any $x \in A$ we have \be \label{eq}wx = \gs(x)w.\ee Suppose that $x$ is an element of $A$ that is not a unit and
set $I = Jx$.  Then  we have a short exact sequence $$0\lra L\lra M\lra N \lra 0,$$  where $L = Ax/I, M = A/I$ and $N = A/Ax$.
\bl $L \cong A/J$ is a simple $A$-module. \el \bpf The map $f$ from $A$ to $L$ sending $a$ to $ax + Jx$ is clearly surjective with kernel containing $J$.  If $a \in \Ker f$, then $(a-j)x = 0$ for some $j\in J,$ whence $a \in J.$ \epf
\noi An interesting feature of our construction is that remaining assumptions involve only $L$ and $N$.
There is a single additional assumption on $L.$
\begin{itemize}
\item[{(4)}] For all $m \ge 0$ the equation
\be \label{e201} \gs^m(x)a - 1 \in J\ee \ei
has no solution for $a\in A.$ For $z \in A$, denote the image of $z$ in $M = A/I$ by ${\overline{z}}$. Then equation (\ref{e201}) is equivalent to
\be \label{e202} \gs^m(x)a{\overline{x}} = \overline{x}\ee and equation (\ref{e202}) always has a solution if $L$ is divisible.  Since $L$ obviously cannot be injective, some condition similar to (4) must be necessary if our construction is to go through.\\ \\
Finally we make the following assumptions on  $N.$
\begin{itemize}
\item[{(5)}] $N$ has a strictly descending chain of submodules
\be \label{sdc}N \supset wN \supset \ldots \supset  w^mN \supset \ldots\ee 
\item[{(6)}] Every nonzero submodule of $N$ contains $w^mN $ for some $m.$
\end{itemize}
\bt \label{mt}
Under assumptions $\mbox{\rm  (1)-(6)}$, $M$ is an essential extension of $L.$\et \bpf  Note that the assumptions are unchanged if we replace $w$ by $\mu^{-1}w$.  Thus we can assume that $\mu = 1.$ Suppose $U$ is a left ideal of $A$ strictly containing $I.$ We need to show that $U$ contains $Ax.$ It follows easily from (6) that  $U$ contains an element of the from $w^m - ax$ for some $a\in A.$ Set $y = \gs^m(x)$. Then from (\ref{eq}) and  (3) we  have
\begin{eqnarray}\label{e101}
y(w^m - ax) &=&  (w^m - 1)x + (1 -  ya)x \nonumber \\ &\equiv & (1 - ya)x \mod Jx.
\end{eqnarray}
For $z \in A$, denote the image of $z$ in $M = A/I$ by ${\overline{z}}$. From (\ref{e101}) and assumption (4) we have $0 \neq (1 - ya){\overline{x}} \in {\overline{U}} \cap A{\overline{x}},$ so as $L = A{\overline{x}}$ is simple it follows that $A{\overline{x}} \subseteq {\overline{U}}.$ The result follows easily.\epf
\section{Examples (B.1)-(B.3).}
To check assumption (4) we use the following easy result.
\bl \label{ca4} If  for all $m \ge 0,$ there is a subring $B$ of $A$ such that $A = B \oplus J$, and $\gs^m(x) \in B$,then assumption (4) holds.\el
\bpf If $\gs^m(x)a - 1 \in J,$ write $a = b+j$ with $b\in B$ and $j \in J.$ Then $\gs^m(x)b - 1 \in J \cap B = 0$, whence $\gs^m(x)$ is a unit in $A$ a contradiction, since $x$ is assumed to be a non-unit.\epf \noi It is not always possible to choose $B$ to be $\gs$-invariant in Lemma \ref{ca4}. From Theorem \ref{mt} and the next two results, we obtain the non-artinian, monolithic modules in \cite{CM2} Theorems 3.1 and 4.2. \\ \\
 Let $A = K[a, b]$  be the coordinate ring of the quantum plane, as in (B.1) where $ab = q ba$ and  $q \in K \backslash \{0\}$ is a not root of unity. Let $w = ab$ and $J = A(ab-1)$, $B = K[a]$ and $x = a-1 \in B$. Then $w$ is a  normal element and the automorphism $\gs$ determined by equation (\ref{eq}) satisfies $\gs(a) = q^{-1}a$ and $\gs(b) = qb$.  
\bp \label{p1}  \item[{(a)}] $J$ is a maximal left ideal of $A$ and assumption (4) holds.\item[{(b)}] If $N = A/Ax,$ then $N$ is non-artinian, and a complete list of non-zero submodules of $N$ is given by equation (\ref{sdc}).\ep
 \bpf 
Since $A = B \oplus J$ and $\gs$ preserves $B$, the result follows from Steps 1 and 2 in the proof  of
\cite{CM2} Theorem 3.1.
 \epf
\noi Let $A = K[a, b]$  be the quantized Weyl algebra, as in (B.2) where $ab - q ba = 1$ and  $q \in K \backslash \{0\}$ is a not root of unity.  If  $w = ab-ba$, then $w$ is a normal element of $A$ and $w-1 = (q-1)ba \in J = Aa$. The automorphism $\gs$ determined by equation (\ref{eq})  satisfies $\gs(a) = q^{-1}a$ and $\gs(b) = qb$. We have $A = B \oplus J$ with $B = K[b],$ and $\gs(B) = B.$   Let $x = (1-q)b-1 \in B$.
\bp \label{p2}  \item[{(a)}] $J$ is a maximal left ideal of $A$ and assumption (4) holds.\item[{(b)}] If $N = A/Ax,$ then $N$ is non-artinian, and a complete list of non-zero submodules of $N$is given by equation (\ref{sdc}).\ep
 \bpf 
By \cite{CM2} Lemma 4.1, $J$ is a maximal left ideal of $A$, and (4) follows as before. Note that $N \cong K[a]$ as a $K[a]$-module.  Let $u_0 = 1 +Ax$,
and define inductively $u_{n+1} = (q^{-n} a - 1)u_n.$  Then $$ au_n = q^n(u_n + u_{n+1})\;\; \mbox{   and   }\; \;  bu_n=\frac{q^{-n}}{1-q}u_n.$$  Thus (b) follows as in the proof of \cite{CM2} Theorem 4.2 (b).
\epf
Next we show that certain Ore extensions with Gelfand-Kirillov dimension 2 do not have property $(\diamond)$. Assume that $K$ 
has characteristic zero, and  let $d$ be the derivation of the polynomial algebra $K[a]$ determined by $d(a) = a^r$ where $r \ge 1.$ Let $A = K[a, b]$  be the resulting Ore extension, where for $p \in K[a],$ \be \label{aba} pb = bp + d(p).\ee  In particular $$ab = ba + a^{r}.$$ Thus if $w = a$, then $w$ is a  normal element and the automorphism $\gs$ determined by equation (\ref{eq}) satisfies $\gs(a) = a$ and $\gs(b) = b + a^{r-1}$.
We show below that $A$ does not have property $(\diamond)$. When $r = 1,$ $A$ is isomorphic to he enveloping algebra $U(\fb)$ where,  $\fb$ is a Borel subalgebra of $\fsl(2,K)$.  Now by \cite{BGR} Lemma 6.12, if $K$ is algebraically closed, then   any finite dimensional solvable Lie algebra which is non-nilpotent has $\fb$ as an image, and thus we recover the result in (B.3).
\bl \label{wkf} Any ideal invariant under $d$
is generated by a power of $a$.\el \bpf This follows from the well known fact that if an ideal $Q$ is invariant under a derivation, then so too are all the prime ideals that are minimal over $Q$, see for example \cite{BGR} Lemma 4.1.\epf
\noi  Let  $J = A(a-1)$ and $x = b-1$.   
\bp \label{p3}  \item[{(a)}] $J$ is a maximal left ideal of $A$ and assumption (4) holds.\item[{(b)}] If $N = A/Ax,$ then $N$ is non-artinian, and a complete list of non-zero submodules of $N$is given by equation (\ref{sdc}).
\item[{(c)}] The submodules of $N$ are pairwise non-isomorphic.
\ep
 \bpf (a) Set $v_n = b^n + J.$ The elements $\{v_n\}_{n \ge 0}$ form a basis for $A/J,$ and  $av_0 = v_0$.  Assume by induction that
\be \label{e1}(a-1)^nv_n = n!v_0.\ee Then by equation (\ref{aba}), we have
\begin{eqnarray}\label{}
(a-1)^{n+1}v_{n+1} &=& (a-1)^{n+1}bv_{n}
\nonumber \\ &=& b[(a-1)^{n+1} +(n+1)a^r(a-1)^{n}]v_n \nonumber \\ &=& (n+1)!v_0
\nonumber
\end{eqnarray}
It follows easily from equation (\ref{e1}) that $A/J$ is simple. Since $\gs^m(x) = b - 1 + ma^{r-1}$ we have  $A = B \oplus J$ where $B = K[\gs^m(x)]$, thus (4) holds. \\ \\(b) 
Since $A = K[a] \oplus Ax$,  we can identify $N$ with $K[a]$ as a $K[a]$-module.  Suppose $N'$ is a submodule of $N$, and $N' = pK[a]$ for some   $p \in K[a].$ Then
\begin{eqnarray}\label{}
bp&=& pb - d(p)
\nonumber \\ &\equiv& p - d(p) \mod Ax,
\nonumber
\end{eqnarray} and hence $d(p) \in pK[a]$.  Thus (b) follows from Lemma \ref{wkf}.\\ \\
(c) As above we identify $N$ with $K[a]$. If $\phi:a^mN \lra a^{m_1}N$ is an isomorphism, then $\phi(a^m) = a^{m_1}q(a)$ for some polynomial $q$ with $q(0) \neq 0.$ Thus 
\begin{eqnarray}
\label{ } 
 \phi(ba^m) &=& \phi(a^m - m a^{m+r-1})
\nonumber \\
&=& (1 - m a^{r-1})a^{m_1}q(a)
.\nonumber
\end{eqnarray}
and 
\begin{eqnarray}
\label{ } 
b\phi(a^m) &=& b(a^{m_1}q(a))\nonumber \\
&=& a^{m_1}q(a)- a^r(a^{m_1}q(a))'
.\nonumber
\end{eqnarray}
This easily gives
$$(m_1 - m)a^{m_1+r-1}q(a) = a^{m_1+r}q'(a).$$  Now we must have $m = m_1$ since otherwise the left side has 0 as a root of multiplicity at most $m -r+1$, whereas the right side has 0 as a root of multiplicity at least $m -r$.
 \epf
\section{Down-up Algebras.} Given a field $K$ and
$\alpha ,\gb,\gamma$ elements of $K$, the associative algebra
$A=A(\ga,\gb,\gamma)$ over $K$ with generators $d,u$ and defining
relations
               $$(R1)\qquad d^2u=\alpha  dud+\gb ud^2+\gamma d$$
               $$(R2)\qquad du^2=\alpha  udu+\gb u^2d+\gamma u$$
is called a down-up algebra. Down-up algebras were introduced by G. Benkart and
T. Roby \cite{BR}, \cite{BR1}.
In \cite{KirkmanMussonPassman} it is shown that $A=A(\alpha ,\gb,\gamma)$ is noetherian if and only if $\gb\neq 0$, and that these conditions are equivalent to $A$ being a domain.
The main result of this section is as follows.
\bt \label{mth} If $A(\alpha, \beta, \gamma)$ is  a noetherian down-up algebra over a field $K$ of characteristic zero,
then any finitely generated monolithic $A(\alpha, \beta, \gamma)$-module is artinian if
and only if the roots of $X^2-\alpha X-\beta$ are roots of unity.\et
From now on we assume that $X^2 -\ga X -\gb = (X-1)(X-\eta)$ where
$\eta =-\gb$ is not a root of 1, and that $\gb \neq 0.$ Thus $A(\alpha, \beta, \gamma)$
is a Noetherian domain by the above remarks, and $\ga +\gb=1.$  In addition we assume that $\gc \neq 0$. Hence $A(\alpha, \beta, \gamma)$ is isomorphic to a down-up algebra
$$A_\eta = A(1+\eta,-\eta,1).$$  
To prove Theorem \ref{mth} it is enough to prove the result below, as noted in 
\cite{CM2}.\bt \label{dup}If $\eta$ is not a root of unity, then $A_\eta$  does not  have property $(\diamond)$.\et  \noi For the remainder of this section we assume that  $A = A_\eta$ as in Theorem \ref{dup}. We begin with some consequences of (R1) and (R2).
Since $\eta \neq 1,$ we have $\ga \neq 2$. Set $\gep  = (\ga-2)^{-1},$ and $\phi = 1 - \ga\gep = -2(\ga -2)^{-1}.$   As noted in \cite{CM1} Section 1.4 Case 2, the element $w = -ud + du +\gep$ satisfies $$dw = \eta wd, \quad uw = \eta^{-1} wu,$$
  and hence $wA = Aw.$
We remark that $A/Aw$ is isomorphic to the first Weyl algebra (this fact is not used below).
\bl For $n \ge 1,$ we have \be \label{ee3} d u^{2n} = u^{2n}d + n\phi u^{2n-1} + \ga \sum_{i=0}^{n-1} \eta^{-2i-1}w  u^{2n-1}   \ee and
for $n \ge 0,$ \be \label{e4} d u^{2n+1} = u^{2n+1}d +
u^{2n}w  +(n\phi  - \gep)u^{2n} + \ga \sum_{i=0}^{n-1} {\eta^{-1-2i}w } u^{2n}.   \ee
\el \bpf We have
\be \label{ee1} du = w + ud -\gep.\ee Using (R2), then (\ref{ee1}) and the fact that $\ga+ \gb = 1$, we see that  for $j \ge 2,$
\begin{eqnarray*}
d u^{j} & = & [\ga udu + \gb u^2d +  u] u^{j-2}\nonumber \\
& = & [\ga u(w + ud -\gep) + \gb u^2d +  u] u^{j-2}\nonumber \\
& = & [(\ga+\gb)u^2d +  \ga uw +(1 -\ga \gep)u] u^{j-2}\nonumber \\
& = &  u^2d u^{j-2} + \ga uwu^{j-2} +\phi u^{j-1}.
\end{eqnarray*} The result follows easily by induction.\epf
Consider the module $N = A/A(d-1),$ and if $a \in A,$ denote the image of $a$ in $N$ by $\overline{a}$.
Then $N$  has a basis $w^i\overline{u}^j$ with $i, j \ge 0.$
Thus if $B = K[u,w],$ then $N \cong B$ as a left $B$-module. Since $dw^m = \eta^m w^md$, 
$N$ has a strictly descending chain of submodules  as in Assumption (5).
Next we define a filtration on $N$ by setting $$N_n = \sum_{i=0}^{n}u^i K[\w] = \sum_{i=0}^{n}K[w]\u^i.$$
It follows from (\ref{ee3}) and (\ref{e4}) that $dN_n \subseteq N_n.$
Also for $f \in K[w],$ we have \be \label{e5} df(w)\u^n \equiv f(\eta w)\u^n \mod N_{n-1}.  \ee
\bl \label{e6} If $U$ is a non-zero submodule of $N$, then $U$ contains $\w^m$ for some $m.$ \el
\bpf  Suppose that $n$ is minimal such that $U \cap N_n \neq 0.$ We claim that $n = 0.$
If this is not the case then $U +N_{n-1}$ contains an element of the form $x = f(w)\u^n$ for some non-zero polynomial $f$.
Write $f(w) = \sum_{i=r}^s a_iw^i,$ where $a_r \neq 0 \neq a_s.$
If $r < s$, then $U +N_{n-1}$  contains an element of the form $y = w^r\u^n,$ because $\prod_{i=r+1}^s (d - \eta^i)x \in  U +N_{n-1}$.
Thus if $n = 2m$ is even, we can assume that
\be \label{e7} y = w^r\u^{2m} + \sum_{i=0}^{2m-1} g_i(w)\u^i \in  U. \nonumber\ee
Then \be \label{e8} (d-\eta^r)y \equiv  [\eta^rw^r(m\phi + \ga \sum_{i=0}^{n-1} \eta^{-2i-1}w) + g_{ 2m-1}(\eta w) -\eta^r g_{ 2m-1}(w) ]\u^{2m-1}\mod N_{n-2}.\nonumber \ee
By the choice of $n, (d-\eta^r)y $ must be zero mod $N_{n-2}$.
Note that the coefficient of $w^r$ in $g_{ 2m-1}(\eta w) -\eta^r g_{ 2m-1}(w) $ is zero. Thus  looking at the coefficient of $w^r\u^{2m-1}$ on the right side 
above yields $m\phi = 0$, which is a contradiction.
Thus   $n = 2m+1$ is odd, and  we can assume that \be \label{e71} y = w^r\u^n + \sum_{i=0}^{2m} f_i(w)\u^i \in  U +N_{n-2}. \nonumber\ee
Then mod $N_{n-2},$ 
\begin{eqnarray*}\label{e81} (d-\eta^r)y &\equiv & \eta^r w^r[u^{2m}\w +(m\phi -\gep)\u^{2m} 
+ \ga \sum_{i=0}^{m-1} {\eta^{-1-2i}w } \u^{2m}]
\nonumber \\
& +& [f_{2m}(\eta w)-\eta^rf_{2m}(w)]\u^{2m} . 
\end{eqnarray*} 
 By the choice of $n, (d-\eta^r)y $ must be zero mod $N_{n-2}$.
Then  looking at the coefficient of $w^r\u^{2m}$ 
we obtain $m\phi = \gep$ which leads to the contradiction $2m + 1 = 0.$ Thus $U$ contains an element of the form $f(\w)$ with $f \neq 0,$ and the result follows easily.\epf
\noi We have verified assumptions (5) and (6) for the module $N$, and we now turn our attention to the simple module $L.$\\ \\
Following \cite{BR} Proposition 2.2, we define 
the Verma module $V(\lambda)$ with highest weight $\lambda \in K$. Let
$\lambda_{-1}=0, \lambda_0=\lambda$ and for each $n > 0$ set, 
\be \label{rec}\lambda_n=\ga\lambda_{n-1}+\gb\lambda_{n-2}+1.\ee The
Verma module $V(\lambda)$ has basis $\{v_n|n\in {\Bbb N}\}$. The
action of $A$ is defined by
$$d.v_0=0, \mbox{ and }
d.v_n=\lambda_{n-1}v_{n-1}, \mbox{for all}\; n\geq 1$$
$$u.v_n=v_{n+1}.$$

In \cite{BR} Proposition 2.4 it is shown that $V(\lambda)$ is
simple if and only if $\lambda_{n} \neq 0$ for all $n \ge 0$.  Furthermore, by \cite{CM1} Lemma 2.5, $\gl_{n-1} = 0$ if and only if 
 \be \label{rec1}\lambda(\eta-1)=-(1-n(\sum_{i=0}^{n}\eta^i)^{-1}).\ee
\bl \label{sim} 
The algebra $A$ has infinitely many pairwise non-isomorphic simple Verma modules.\el
\bpf The result is evident if $K$ is uncountable, because then we simply require that the highest weight $\gl$ does not satisfy the condition in (\ref{rec1}) for any $n.$  In general we argue as follows.  By \cite{CM1} Proposition 5.5, any Verma module has length at most 3, so by \cite{BR} Proposition 2.23, any Verma module has a simple Verma submodule. Also if $V(\gl)$ is not simple this submodule is generated by $v_n$ where $n$ is the largest integer  such that $\gl_{n-1} = 0$.  This submodule is isomorphic to $V(\gl_n)$. Note that the case covered by \cite{BR1} 
does not arise here.  Now if $\mu = \gl_{n}$ and  $V(\mu)$  is simple, we can solve the recurrence (\ref{rec}) in reverse to find all Verma modules  $V(\gl)$ containing as a $V(\mu)$  simple submodule. Since there can be at most 3 such $\gl$ and $K$ is infinite, the result follows.\epf
\noi Unfortunately it does not seem possible to verify assumption (4) for a simple Verma module. Instead we consider the universal lowest weight modules
$W(\gk)$ defined in \cite{BR} Proposition 2.30 (a).  For $\gk \in K,$   set $\gk_{-1}=0, \gk_0=\gk$ and define for each $n>0$, 
\be \label{rec2}\gk_n= \eta^{-1}(\ga\gk_{n-1}-\gk_{n-2}+1).\ee 
Then $W(\gk)$ has basis $\{a_n|n\in {\Bbb N}\}$. The
action of $A$ is defined as follows, 
$$u.a_0=0, \mbox{ and }
u.a_n=\gk_{n-1}a_{n-1}, \mbox{for all}\; n\geq 1$$
$$d.a_n=a_{n+1}.$$
\bc \label{cor} The algebra $A$ has infinitely many pairwise non-isomorphic simple lowest weight modules
$W(\gk)$. \ec
\bpf By \cite{CM1} Lemma 4.1, there is an isomorphism from $A$ onto $A' = A_{\eta^{-1}}$ which interchanges the generators $u$ and $d.$ Under this isomorphism, any
Verma module for $A'$ becomes a module of the form $W(\gk)$ for $A$, so the result follows.\epf
\noi {\it Proof of Theorem \ref{dup}.} Let $L=W(\gk)$  be a simple lowest weight module, and let $J$ be the annihilator of the lowest weight vector $a_0$ in $A$.  Then $J = Au + A(ud -\gk).$ The normal element $w = -ud + du +\gep$ satisfies $w -\mu \in J$ where $\mu = -\gk+\gep.$ By Corollary \ref{cor} we can arrange that $\mu$ is non-zero.  Set $x = d-1.$ It only remains to check assumption (4).  This holds because $A = B \oplus J$ with $x \in B =K[d],$ and $B$ is $\gs$-invariant.\hfill $\Box$   
\section{Remarks and Problems.}

\begin{itemize}
 \item[{(a)}]  We call  a finitely generated module $E$ over a left noetherian ring {\it uniserial} if the submodules of $E$ are totally ordered by inclusion.  For $E$ uniserial define a descending chain of submodules $\{E_\ga\}$ as follows.  For any ordinal $\ga,$ if $E_\ga \neq 0$ let $E_{\ga + 1}$ be the unique maximal submodule  of $E_\ga$.  For a limit ordinal $\gb$ such that $E_\ga \neq 0$  for $\ga < \gb,$ set $E_\gb = \bigcap_{\ga < \gb} E_\ga.$ There is a smallest ordinal $\gt$ such that $E_\gt = 0$, and we call  $\gt$ the {\it depth} of $E$. As noted in the introduction to  \cite{J2}, it follows from \cite{J2} Theorem 4.6, that for any ordinal $\gt$ there is a left noetherian ring $A$ such that the left regular module is uniserial with depth $\gt.$  The modules  $M$ constructed using Theorem \ref{mt} with the aid of the results in Section 3 are all uniserial with 
depth $\go+1$ where $\go$ is the first infinite ordinal. What other module depths are possible for uniserial modules over (two-sided) noetherian rings?
\item[{(b)}] 
If $N$ is as in Propositions \ref{p1} and \ref{p2} (resp. \ref{p3}), then $N$ is incompressible and critical by \cite{CM2} Theorems 3.1 and 4.2, (resp. Proposition \ref{p3} (c)).   The first example of an  a incompressible and critical module was found by Ken Goodearl,  see \cite{G}, to which we refer for the definitions. 
In general is there a connection between rings that do not have $A$ property $(\diamond)$, and incompressible critical modules? \item[{(c)}] Suppose that  $A$ is  a  Noetherian ring, and $P$ an ideal of $A$ such that $A/P$ is simple artinian with simple module $S$. Is  the injective hull of a $S$ as an $A$-module locally artinian?
\item[{(d)}] Define a noetherian ring $A$ to be $(\diamond)$  {\it extremal} if it does not have property $(\diamond),$  but every proper homomorphic image has property $(\diamond).$ What can be said about $(\diamond)$  extremal rings?  If $A$ is an algebra over a field having finite Gelfand-Kirillov dimension and $A$ is $(\diamond)$  extremal, must $A$ be prime?  The Goodearl-Schofield  example shows that this is not true without the GK dimension  hypothesis.  It seems likely that the algebra $A_\eta$ in Theorem \ref{dup} is $(\diamond)$  extremal. We note the following result.\ei
\bp \label{er} Suppose that $A$ is a $K$-algebra such that the endomorphism ring of every simple $A$-module is algebraic over $K.$ If $A$ is $(\diamond)$  extremal the center $Z$ of $A$ is algebraic over $K.$  \ep \bpf If $Z$ is not algebraic over $K$ then, for every simple module $L$, the natural map $Z \lra \End_A L$ has non-zero kernel {\bf m}. Then if the injective hull of $L$ as an $A/{\bf m}A$ is locally artinian, then so too is its injective hull over $A$, see \cite{CarvalhoLompPusat} Proposition 1.6. Thus $A$ cannot be $(\diamond)$  extremal.\epf The hypothesis that the endomorphism ring of every simple $A$-module is algebraic over $K$ is known to hold for many algebras, for example it holds for almost commutative algebras (Quillen's Lemma) and for an algebra of countable dimension over an uncountable field. 
\pagebreak
 
\label{refs}

\end{document}